\documentclass[11pt,A4,twoside]{article}
\usepackage[latin1]{inputenc}
\usepackage{amssymb}
\usepackage{amsmath}
\usepackage{amsfonts}
\usepackage{amsbsy}
\usepackage{natbib}
\usepackage{eucal}
\usepackage{graphicx}
\usepackage{titlesec}
\usepackage[english]{babel}
\usepackage{times}
\usepackage{titling}
\usepackage{mathptmx}
\usepackage{float}
\usepackage[font=small]{caption}
\usepackage{array}
\usepackage[colorlinks,citecolor=blue,urlcolor=blue]{hyperref}
\usepackage{authblk}
\usepackage{amsthm}
\usepackage{pdfpages}
\usepackage{url}
\usepackage{setspace}
\DeclareSymbolFont{txgreek}{OML}{cmr}{m}{it}

\renewcommand{\abstract}[1]{{\small\noindent
\hrulefill\par \vspace*{0.1cm}\noindent{\small\bf\sffamily
{Abstract}}\parindent=0pt\par\noindent\vspace{-0.1cm}\noindent\hrulefill\par\vspace*{0.5\baselineskip}\hspace*{0cm}\renewcommand{\baselinestretch}{1.1}\sffamily{#1}\par\vspace*{-0.1cm}\noindent\hrulefill}}

\newtheorem{theorem}{Theorem}[section]
\newtheorem{corollary}{Corollary}[theorem]

\def\and{,\;}
\DeclareMathSizes{11}{11}{7.8}{6.5}

\def\paragraf{\fontsize{9}{10pt}\fontfamily{phv}\fontshape{it}\selectfont}
\titleformat{\paragraph}
{\titlerule[0pt]
\vspace{0cm}%
\paragraf}
{\theparagraph}{.5em}{\vspace*{0\baselineskip}}

\def\titol{\fontsize{12.045}{12pt}\fontfamily{phv}\fontseries{b}\selectfont}
\titleformat{\section}
{\titlerule[0pt]
\vspace{0.2cm}%
\titol}
{\thesection.}{.5em}{\vspace*{0\baselineskip}}

\def\titolp{\fontsize{11.045}{11pt}\fontfamily{phv}\fontseries{b}\fontshape{it}\selectfont}
\titleformat{\subsection}
{\titlerule[0pt]\bigskip
\vspace{-0.4cm}%
\titolp}
{\thesubsection.}{.5em}{\vspace*{0\baselineskip}}

\def\titolpp{\fontsize{10.045}{10pt}\fontfamily{phv}\fontshape{it}\selectfont}
\titleformat{\subsubsection}
{\titlerule[0pt]
\vspace{-0.4cm}%
\titolpp}
{\thesubsubsection.}{.5em}{\vspace*{0\baselineskip}}

\usepackage{titling}
    \pretitle{\begin{center}\sffamily\fontsize{18pt}{20pt}\selectfont}
    \posttitle{\par\end{center}}
    \preauthor{\begin{center}\fontsize{12pt}{14pt}\selectfont}
    \postauthor{\par\end{center}\vspace{0bp}}
    \predate{}
    \date{}
    \postdate{}

\title{Application of  Random Walk in Manpower Planning}

\thanksmarkseries{arabic}
\author{Theo van Uem\thanks{Amsterdam University of Applied Sciences, Amsterdam, The Netherlands. tjvanuem@gmail.com} }

\def\headers#1{\fontsize{8.5}{10}\centering\sffamily\itshape{#1}}
\def\page#1{\fontsize{8.5}{10}\sffamily{#1}}
\usepackage{fancyhdr}%Headers

\begin{document}

\maketitle

% Headers
\thispagestyle{empty}
\renewcommand{\headrulewidth}{0truecm}
\pagestyle{fancy}
\rhead[\headers{Application of  Random Walk in Manpower Planning}]{\page{\thepage}}
\lhead[\page{\thepage}]{\headers{Theo van Uem}}
 \lfoot{} \rfoot{}
\cfoot{}

\abstract{The career of an employee can be described (under certain circumstances) by a random walk, where the states of the random walk are determined by the level and position of an employee. At each decision moment the state of the employee is changed by four stochastic transformations:  upgrading one position at the same level, upgrading one level, staying until the next decision moment in the current state and absorption in the current state. We obtain explicit formula for the long term behavior of the distribution of all employees using generating functions.    }


\renewcommand{\baselinestretch}{1.2}
\bigskip

%\hbadness=99999

 % Keywords
                       % Math. Subj. Class. codes

\section{Introduction}

By The National Institutes of Health (https://hr.nih.gov/workforce/workforce-planning) is Workforce Planning the process of analyzing, forecasting, and planning workforce supply and demand, assessing gaps, and determining target talent management interventions to ensure that an organization has the right people - with the right skills in the right places at the right time - to fulfill its mandate and strategic objectives.\\
In this paper we focus on forecasting and planning.\\
Markov chains can be used to forecast probability distributions for different cohorts in an organization. (see e.g.\citep{9} and \citep{10}).
In this paper we suppose that the total number of employees in an organization is more or less constant. Employees who leave the organization are replaced by employees of the same level. The situation within such an organization is not static: employees can grow in their own level or  jump to a higher level.
The model we use is a random walk model where the states of the system are defined by level and position (within that level) of the employee. Let $t_k  \  ( k \in \mathbb {N})$ be the times where decisions about changing level and/or positions are made.
Both levels and positions are unlimited in our model, but in practice we get always finite solutions. Let an employee be in state $(n,m)$  (m is level and n is the position within that level) at time $t_k$. Then at time $t_{k+1}$ our employee will:
\begin{itemize}
\item jump to state $(n+1,m$) with probability $p$ (upgrading one position at the same level)
\item jump to state $(0,m+1)$ with probability $r$ (upgrading one level, starting in position 0)
\item stay until $t_{k+1}$ in $(n,m)$ with probability $s$ 
\item absorb in $(n,m)$ with probability $t$ 
\end{itemize}
where $p+r+s+t=1$ and we demand $prst>0$ \\
By choosing values of $p,r,s,t$ we can see, using our results,  what the long term effect is and we can implement desired values of $p,r,s,t$  in the practice of our organization.

\section{Generating functions}
In general, for a discrete Markov chain,  we define the expected number of visits to state j when starting in state i by:
\begin{equation*}
x_j=x_{i,j}=\sum_{k=0}^{\infty}p_{i,j}^{(k)}
\end{equation*}
where $ p_{i,j}^{(k)}$ is the probability of moving from state $i$ to state $j$ in $k$ steps.\\
In our situation we start at level $0$  in situation  $i_0.$\\
Let $f_{n,m}$ be the expected number of visits to situation $n$ at level $m$ ($n\in \mathbb {N}, m\in \mathbb{N})$.\\
Last step analysis gives:\\
 For level $m$ we have:
\begin{equation*}
f_{k,0}=0   \quad  ( k=0,1..,i_0-1) \quad (m=0)
\end{equation*}
\begin{equation}
\label{1}
f_{i_0,0}=sf_{i_0,0}+1   \quad  ( m=0) 
\end{equation}
\begin{equation}
\label{2}
f_{0,m}=r\sum_{k=0}^{\infty}f_{k,m-1}+sf_{0,m} \quad ( m \geq 1)\\
\end{equation}
\begin{equation}
\label{3}
f_{n,m}=pf_{n-1,m}+sf_{n,m} \quad  ( m \geq 1, \ n \geq 1) 
\end{equation}

 We define generating functions
\begin{equation*}
F(u,v)=\sum_{n=0}^{\infty}\sum_{m=0}^{\infty}f_{n,m}u^nv^m \quad (0<u\leq 1, 0<v\leq 1)
\end{equation*}
\begin{equation*}
G_{m}(u)=\sum_{n=0}^{\infty}f_{n,m}u^n\quad (0<u\leq 1, m\in \mathbb {N})
\end{equation*}
Using \eqref{1} we get:
\begin{equation}
f_{i_0,0}=\frac{1}{1-s}
\end{equation}
\begin{theorem}
\begin{equation}
\label{vijf}
G_0(u)=\frac{u^{i_0}}{(1-s-pu)}
\end{equation}
\begin{equation}
\label{5}
G_m(u)=\frac{1}{(1-s-pu)}\left (\frac{r}{r+t}\right )^{m} \quad (m=1,2,..)
\end{equation}
\end{theorem}
\begin{proof}
CASE $m=0$.\\
Using \eqref{3} and \eqref{1} we get:
\begin{equation*}
\begin{split}
G_0(u)&=\sum_{n=0}^{\infty}f_{n,0}u^n\\
&=\sum_{n=i_0}^{\infty}f_{n,0}u^n\\
&=f_{i_0,0}u^{i_0}+\sum_{n=i_0+1}^{\infty}\{pf_{n-1,0}+sf_{n,0}\}u^n\\
&=f_{i_0,0}u^{i_0}+pu\sum_{n=i_0+1}^{\infty}f_{n-1,0}u^{n-1}+s\sum_{n=i_0+1}^{\infty}f_{n,0}u^n\\
&=f_{i_0,0}u^{i_0}+puG_{0}(u)+s\{G_{0}(u)-f_{i_0,0}u^{i_0}\}\\
&=u^{i_0}+(pu+s)G_{0}(u)
\end{split}
\end{equation*}
  
CASE $m>0$.\\
Using \eqref{3}  we get:
\begin{equation*}
\begin{split}
G_{m}(u)&=\sum_{n=0}^{\infty}f_{n,m}u^n\\
&=f_{0,m}+\sum_{n=1}^{\infty}f_{n,m}u^n\\
&=f_{0,m}+\frac{p}{1-s}\sum_{n=1}^{\infty}f_{n-1,m}u^n\\
&=f_{0,m}+\frac{pu}{1-s}G_{m}(u)
\end{split}
\end{equation*}
 so:
\begin{equation}
\label{6}
G_m(u)=\frac{f_{0,m}}{1-\frac{pu}{1-s}} \quad (m\geq 1)
\end{equation}
Using \eqref{2} and \eqref{6} gives:\\
$
f_{0,m}=\frac{r}{1-s}\sum_{k=0}^{\infty}f_{k,m-1}=\frac{r}{1-s}G_{m-1}(1)=\frac{r}{1-s}\frac{f_{0,m-1}}{1-\frac{p}{1-s}}=\frac{rf_{0,m-1}}{1-s-p}=\frac{rf_{0,m-1}}{r+t} \  (m>1)
$\\
and 
$
f_{0,1}=\frac{r}{1-s}\sum_{k=i_0}^{\infty}f_{k,m-1}=\frac{r}{1-s}G_{0}(1)=\frac{r}{r+p}\frac{1}{1-s}
$, so:
\begin{equation}
\label{7}
f_{0,m}=\frac{1}{1-s}(\frac{r}{r+t})^m
\end{equation}
Combining \eqref{6} and \eqref{7} we get \eqref{5}.
\end{proof}
\begin{theorem} 
\begin{equation}
F(u,v)=\frac{(r+t)u^{i_0}+r(1-u^{i_0})v}{(r+t-rv)(1-s-pu)}
\end{equation}
\end{theorem}
\smallskip
\begin{proof}
Using \eqref{3}, \eqref{5} and \eqref{7} we get:
\begin{equation*}
\begin{split}
F(u,v)&=\sum_{n=0}^{\infty}\sum_{m=0}^{\infty}f_{n,m}u^nv^m\\
&=\sum_{n=0}^{\infty}f_{n,0}u^n+\sum_{m=1}^{\infty}f_{0,m}v^m+\sum_{n=1}^{\infty}\sum_{m=1}^{\infty}f_{n,m}u^nv^m\\
&=G_0(u)+\sum_{m=1}^{\infty} \frac{1}{1-s}(\frac{r}{r+t})^mv^m+\sum_{n=1}^{\infty}\sum_{m=1}^{\infty}\frac{p}{1-s}f_{n-1,m}u^nv^m\\
&=G_0(u)+\frac{1}{1-s}\frac{\frac{rv}{r+t}}{1-\frac{rv}{r+t}}+\frac{pu}{1-s}\sum_{n=1}^{\infty}\sum_{m=1}^{\infty}f_{n-1,m}u^{n-1}v^m\\
&=G_0(u)+\frac{rv}{(1-s)(r+t-rv)}+\frac{pu}{1-s}\{F(u,v)-G_0(u)\}
\end{split}
\end{equation*}
so:
$$
(1-\frac{pu}{1-s}) F(u,v)=(1-\frac{pu}{1-s})G_0(u)+\frac{rv}{(1-s)(r+t-rv)}
$$
and:
$$
F(u,v)=\frac{u^{i_0}}{1-s-pu}+\frac{rv}{(r+t-rv)(1-s-pu)}=\frac{(r+t)u^{i_0}+r(1-u^{i_0})v}{(r+t-rv)(1-s-pu)}
$$
\end{proof}

\section{Absorption probabilities}
Let $\theta=\frac{r}{r+t}$ and $\omega=\frac{p}{1-s}$.
The absorption probability in state $(n,m)$, where $m$ is the level and $n$ the position, is $tf_{n,m}$, where:
\begin{theorem} 
\begin{equation}
f_{n,0}=\frac{\omega^{n-i_0}}{1-s} \quad (n=i_0,i_0+1,...)
\end{equation}
\begin{equation}
f_{n,m}=\frac{1}{1-s}\theta^m\omega^n  \quad (n\in \mathbb {N},  \ \ m=1,2,..)
\end{equation}
\end{theorem}
\begin{proof}
\begin{equation*}
\begin{split}
F(u,v)&=\sum_{n=0}^{\infty}\sum_{m=0}^{\infty}f_{n,m}u^nv^m\\
&=\frac{(r+t)u^{i_0}+r(1-u^{i_0})v}{(r+t-rv)(1-s-pu)}\\
&=\frac{u^{i_0}+(1-u^{i_0})\theta v}{(1-\theta v)(1-s)(1-\omega u)}
\end{split}
\end{equation*}
so:
\begin{equation*}
\begin{split}
(1-s)F(u,v)&=\{u^{i_0}+(1-u^{i_0})\theta v\} \sum_{n} \omega^n u^n \sum_{m} \theta^m v^m\\
&=\omega ^{-i_0}\sum_n \sum_m \omega ^{n+i_0} u^{n+i_0} \theta ^m v^m+ \sum_n \sum_m \omega ^{n} u^{n} \theta ^{m+1} v^{m+1}\\
&=\omega ^{-i_0} \sum_{n \geq i_0}\sum_{m \geq 0}  \omega ^{n} u^{n} \theta ^{m} v^{m}+ \sum_{n \geq 0}\sum_{m \geq 1}  \omega ^{n} u^{n} \theta ^{m} v^{m}-\omega ^{-i_0} \sum_{n \geq i_0}\sum_{m \geq 1}  \omega ^{n} u^{n} \theta ^{m} v^{m}\\
&= \sum_{n \geq 0}\sum_{m \geq 1}  \omega ^{n} u^{n} \theta ^{m} v^{m}+\omega ^{-i_0} \sum_{n \geq i_0}  \omega ^{n} u^{n} 
\end{split}
\end{equation*}
\end{proof}
\smallskip

The probability of absorption in level $m \ $ is $\pi_m=t\sum_{n}f_{n,m} \  \  (m\in \mathbb{N})$.
\begin{theorem}
\begin{equation}
\pi_m=(1-\theta) \theta ^m \quad (m\in \mathbb {N})
\end{equation}
\end{theorem}
\begin{proof}
\begin{equation*}
\pi_0=\sum_{n=i_0}^{\infty}\frac{t\omega^{n-i_0}}{1-s}=\frac{t}{(1-s)(1-w)}=1-\theta
\end{equation*}
\begin{equation*}
\pi_m=\sum_{n=0}^{\infty}\frac{t\omega^{n}\theta^m}{1-s}=\frac{t\theta^m}{(1-s)(1-w)}=(1-\theta)\theta^m  \quad (m>0)
\end{equation*}
\end{proof}

\begin{corollary}
$\pi_m $ is a (geometric) probability distribution, so absorption always occurs: $\sum_{n=0}^{\infty}\sum_{m=0}^{\infty}tf_{n,m}=1$.
\end{corollary}

\section{Completing the model}
In the preceding sections we constructed the base of our model. We will complete the model. First we define a distribution $d_{n,m}$ as the number of persons in the organization which are on level $m$ and in position $n$. Let $d_{n,m}^{stable}$ be the long run situation (stable distribution). We define $f_{n,m}^{n_0,m_0}$ as the value of $f_{n,m}$ when starting in  $({n_0,m_0})$.

\begin{theorem}
\begin{equation}
d_{n,m}^{stable}=\sum_{n_0}\sum_{m_0}tf_{n,m}^{n_0,m_0}d_{n_0,m_0}
\end{equation}
where:
\begin{equation}
f_{n,m_0}^{n_0,m_0}=\frac{\omega^{n-n_0}}{1-s} \quad (n=n_0,n_0+1,...)
\end{equation}
\begin{equation}
f_{n,m}^{n_0,m_0}=\frac{1}{1-s}\theta^{m-m_0}\omega^n  \quad (n\in \mathbb {N},  \ \ m=m_0+1,m_0+2,..)
\end{equation}
\end{theorem}
\begin{proof}
The results in the preceding sections were obtained for situations where we started in position $i_0$ at level $m=0$, but results are easily adapted when starting in  position $n_0$ at level $m_0$.
\end{proof}

\section{Example}
We give a simple example for an organization with 260 employees. Running the Excel/VBA application RWMP (Random Walk and  Manpower Planning), see appendix \ref{A}, we see the result in appendix \ref{B}, where the first table is the start distribution of the organization and the second table is the long run (stable) distribution.
\appendix
\includepdf[
    pages=1,
    scale=0.7,
    nup=1x1,
    frame,
    offset={2.5cm, -1.0cm},
    pagecommand={%
       \section{VBA-Excel program RWMP}
	\label{A}
}       
]{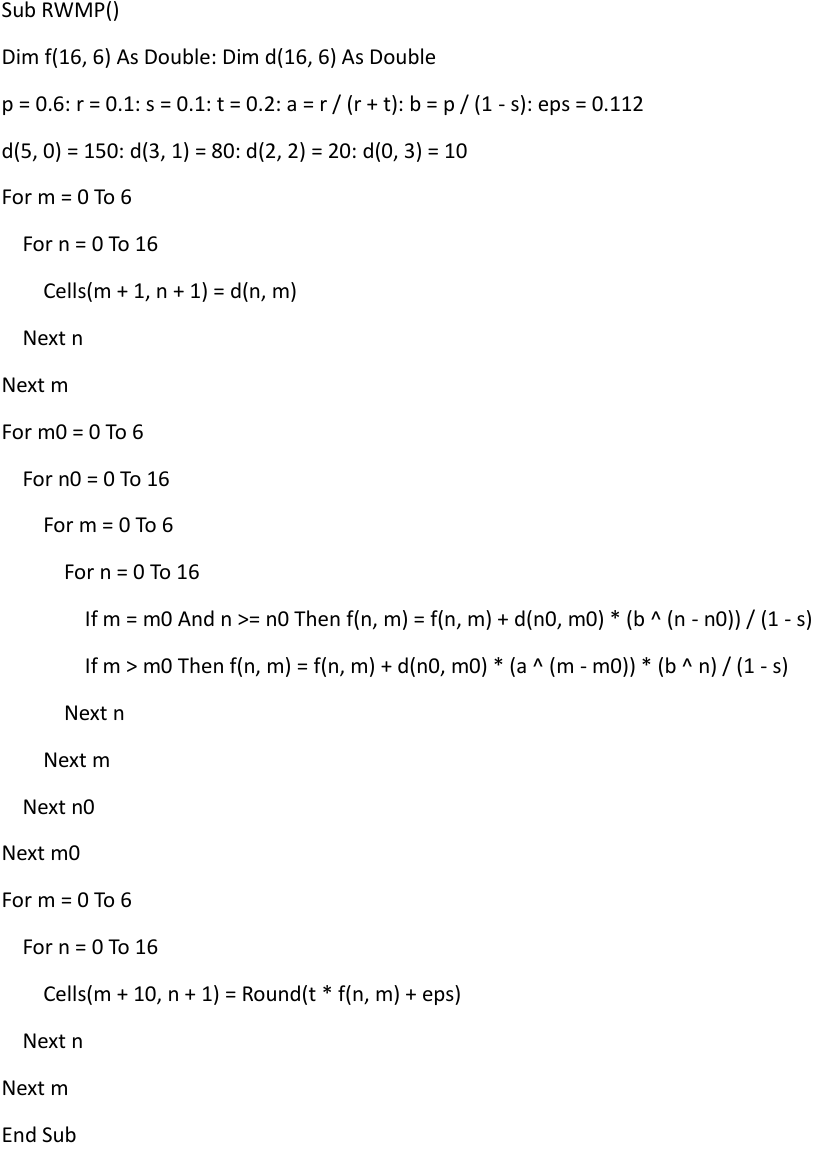}
\section{Output RWMP}
\label{B}
\includegraphics[width=\textwidth]{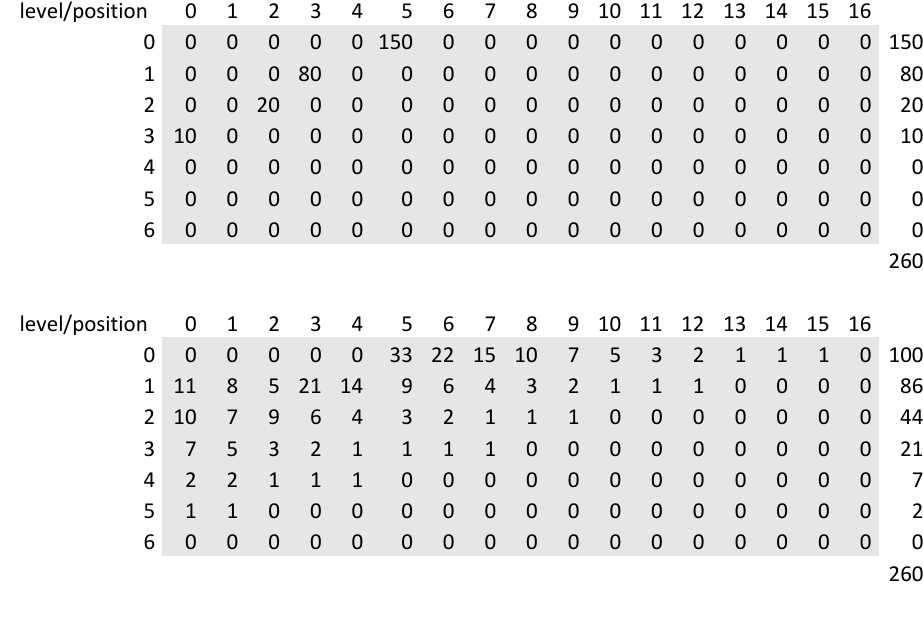}
\clearpage
\bibliographystyle{chicago}

\end{document}